\documentclass{article}
\usepackage{arxiv}
\usepackage{amsthm} 
\usepackage{amssymb}
\usepackage{amsmath}
\usepackage{empheq}
\usepackage{amscd}
\usepackage{graphicx}
\usepackage{graphics}
\usepackage[noadjust]{cite}
\usepackage{amsthm}
\usepackage{caption}
\usepackage{multirow}
\usepackage{array}
\usepackage{rotating}
\usepackage[norelsize,ruled]{algorithm2e}
\usepackage{hyperref}

\usepackage{indentfirst}
\usepackage{tikz}
\usepackage{calc}
\usepackage{epsfig}
\usepackage[numbers]{natbib}
\usepackage[makeroom]{cancel}
\usepackage{multicol}
\usepackage{csquotes}
\usepackage{enumitem}
\usepackage{nicefrac}

\usepackage{hyperref}       
\hypersetup{colorlinks,linkcolor={green},citecolor={green},urlcolor={black}}
\usepackage{optidef}

\newcommand{\bolds}[1]{{#1}}   
\renewcommand{\bold}[1]{{#1}}

\newcommand{\tilt}{\theta}   

\newcommand{\tcm}{}

\renewcommand{\eqref}[1]{equation \ref{#1}}                  

\DeclareMathOperator*{\minimize}{minimize}

\begin{document}
\title{Optimal Space-Variant\\ Anisotropic Tikhonov Regularization\\ for Full Waveform Inversion of Sparse Data}

\author{\href{https://orcid.org/0000-0002-9879-2944}{\hspace{1mm}Ali Gholami} \\
  Institute of Geophysics, Polish Academy of Sciences, Warsaw, Poland\\
  \texttt{agholami@igf.edu.pl} \\ 
  \And 
  \href{https://orcid.org/0000-0002-9879-2944}{\hspace{1mm}Silvia Gazzola}\\
   Silvia Gazzola is with the Department of Mathematics, University of Pisa, Pisa, Italy\\
  \texttt{silvia.gazzola@unipi.it} \\
  }

\graphicspath{{./figs/}}

\renewcommand{\shorttitle}{Optimal Space-Variant Anisotropic Tikhonov Regularization ~~~~~~~~~~~~~ Gholami and Gazzola}

\maketitle

\begin{abstract}
Full waveform inversion (FWI) is a challenging, ill-posed nonlinear inverse problem that requires robust regularization techniques to stabilize the solution and yield geologically meaningful results, especially when dealing with sparse data. Standard Tikhonov regularization, though commonly employed in FWI, applies uniform smoothing that often leads to oversmoothing of key geological features, as it fails to account for the underlying structural complexity of the subsurface. To overcome this limitation, we propose an FWI algorithm enhanced by a novel Tikhonov regularization technique involving a parametric regularizer, which is automatically optimized to apply directional space-variant smoothing. Specifically, the parameters defining the regularizer (orientation and anisotropy) are treated as additional unknowns in the objective function, allowing the algorithm to estimate them simultaneously with the model. We introduce an efficient numerical implementation for FWI with the proposed space-variant regularization. Numerical tests on sparse data demonstrate the proposed method's effectiveness and robustness in reconstructing models with complex structures, significantly improving the inversion results compared to the standard Tikhonov regularization.
\end{abstract}



\graphicspath{{"./figures/"}}
\section{Introduction}
Full waveform inversion (FWI) is a widely used seismic imaging technique with broad applications across various fields of geoscience \cite{Aghamiry_2019_CRO,Jin_2022_RLM,Operto_2023_FWI,Tian_2024_ESW}. FWI enables the recovery of key subsurface physical properties, including velocity, density, attenuation, conductivity, permittivity, and anisotropic parameters. Although the theory behind FWI has been explored for over 30 years, recent advancements in computational resources and numerical optimization techniques have significantly accelerated research and development in this field \cite{Metivier_2022_NAD,Operto_2023_FWI,Barnier_2023_FWIp,Tian_2024_ESW}.

Mathematically, FWI is a large-scale, nonlinear, nonconvex, and ill-posed inverse problem, requiring efficient and robust numerical algorithms for its solution. Given the problem's large scale, gradient-based optimization algorithms are a natural choice for addressing its nonlinearity through iterative linearization. A well-known challenge in FWI is the nonconvexity of the optimization problem, commonly referred to as the ``cycle skipping" issue. Indeed, the optimization landscape of the $\ell_2$ misfit function can be highly nonconvex, containing numerous local minima that can trap local optimization methods. Consequently, the convergence of these algorithms heavily depends on the accuracy of the initial model, especially when the data lacks sufficiently low-frequency content \cite{Operto_2023_FWI,Wang_2020_ELW,Wang_2019_RLW}.
This nonconvexity challenge becomes even more severe in the case of sparse data acquisition or complex geological structures, where inadequate subsurface sampling further exacerbates the problem \cite{Gorszczyk_2017_TRW}.

A common strategy for mitigating the cycle skipping issue in FWI is the use of regularization. This involves adding an additional term to the objective function, transforming the original problem into one that is more well-posed, for example, by making it smoother, more convex, or improving its condition number. Regularization is not a new concept in FWI, and various regularization techniques have been explored to enhance the quality of the reconstructed models (see \cite{Aghamiry_2019_CRO} and references therein).

The most widely adopted regularization method is Tikhonov regularization \cite{Tikhonov_1963_SIF}, which promotes model smoothness by penalizing large gradients in the model. This can be implemented either by explicitly adding the $\ell_2$-norm of the model gradient to the objective function or by smoothing the gradient of the objective function through (convolution) filtering techniques. However, such stationary reguarization leads to over-smoothing of the structured features in the solution. 
Recently, more structure oriented smoothing based on anisotropic diffusion filters has been used for FWI \cite{Metivier_2022_NAD}.
Essentially, by employing nonlinear diffusion filters to smooth the model gradient along the dominant orientations, model fluctuations along the structure are penalized while preserving variations across it, thereby retaining important geological features in the reconstructed model.
Recently,  \cite{Gholami_2024_RES} proposed an anisotropic Tikhonov regularization for linear inverse problems where the local orientation angles are adaptively estimated and updated with the model parameters while the isotropy of the smoothing kernels remains fixed over the space. 
 
In this paper, we introduce an optimal space-variant anisotropic Tikhonov regularization technique, where all the parameters of the smoothing kernels are spatially adaptive and optimized for each pixel of the model. In this innovative approach, local orientation angles and weight parameters are incorporated as additional variables within the FWI objective function. Specifically, the orientation angles determine the preferred structural directions, guiding the penalization of model gradients along these orientations. Meanwhile, the weight parameters define the shape of the smoothing kernels by appropriately scaling the model gradients along the principal directions. 
The objective function is jointly minimized with respect to the orientation field, weights, and the model parameters. Numerical experiments demonstrate that this novel regularization method significantly improves the quality of reconstructed models, especially when compared to conventional Tikhonov regularization or filtering techniques. This enhancement is particularly evident when dealing with sparse data, where the proposed approach effectively preserves complex geological structures.

\section{Theory}
In this section, we propose a FWI method enhanced by space-variant anisotropic Tikhonov regularization. 
\subsection{Full Waveform Inversion}
In the frequency domain, FWI can be formulated as the following nonlinearly constrained optimization problem  \cite{Pratt_1998_GNF,Operto_2023_FWI}:
\begin{equation}\label{main}
    \minimize_{m,u_s}\frac{1}{2}\sum_{s=1}^{n_s}\|Pu_s-d_s\|_2^2~~\text{subject to}~~ A(m)u_s=b_s,
\end{equation}
where $\|\cdot\|_2$ denotes the vector $2$-norm, $\bold{m} \in {\mathbb{R}}^{n \times 1}$ represents the model parameters (the reciprocal of the squared velocity), $\bold{u}_s \in {\mathbb{C}}^{n\times 1}, s=1,...,n_s$,  denotes the wavefields for each source $b_s, s=1,...,n_s$, $\bold{d}_s \in {\mathbb{C}}^{n_r \times 1}$ are the observed data ($n_r$ being the number of receivers), and $\bold{P} \in {\mathbb{R}}^{n_r\times n}$ (with $n_r\ll n$) is the sampling operator.
Each monochromatic wavfield $u_s$ satisfies the discretized wave equation $A(m)u_s=b_s$, where $A(m)=(\omega^2 \text{diag}(m)+\nabla^2) \in \mathbb{C}^{n \times n}$ is the Helmholtz operator. This operator is constructed with sufficient accuracy and appropriate boundary conditions; $\nabla^2$ is the Laplacian operator, and $\omega$ denotes the angular frequency. However, problem \eqref{main} is inherently ill-posed, and robust regularization techniques are essential to ensure convergence to a geologically meaningful solution, particularly when inverting data obtained by sparse source/receivers and starting from an inaccurate initial model.

\subsection{Space-Variant Anisotropic Tikhonov Regularization}
A standard and widely used technique for stabilization of inverse problems is Tikhonov regularization \cite{Tikhonov_1963_SIF}, consisting in modifying the objective function in \eqref{main} by adding a term of the form: 
\begin{equation} \label{Reg0}
\mathcal{R}(m)=\frac12\|\nabla \bold{m}\|_2^2=\frac12\sum_{i=1}^N\| [\nabla \bold{m}]_i\|_2^2,
\end{equation}
where
\begin{equation}\label{grad_comp}
[\nabla\bold{m}]_i=
\begin{bmatrix}
[\nabla_{\!x}\bold{m}]_i\\
[\nabla_{\!z}\bold{m}]_i
\end{bmatrix},
\end{equation}
 $\nabla_{\!x}$ and $\nabla_{\!z}$ are finite difference operators discretizing the first-order horizontal ($x$) and vertical ($z$) derivatives. 
Tikhonov regularization discourages large gradients in the solution. However, due to its sum-of-squares structure, it is most effective when the model's gradient field follows a Gaussian distribution \cite{Tarantola_2005_IPT}. For more complex models, such as those encountered in geophysics, this assumption often does not hold. As a result, Tikhonov regularization tends to penalize larger gradient components too severely, leading to oversmoothing and suboptimal performance. Therefore, it is essential to optimize and adapt the regularization approach to account for the structural complexity of the subsurface.

 \begin{figure}
 \center
 \includegraphics[width=0.70\columnwidth]{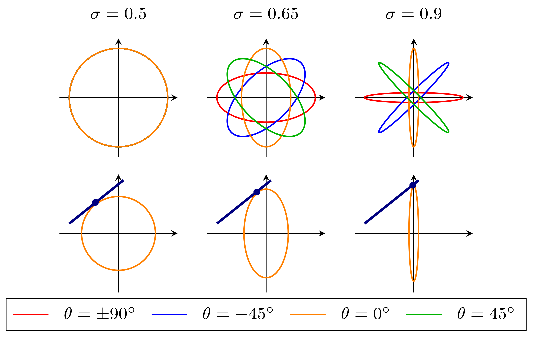}
 \caption{
\tcm{Top row: Geometric illustration of the effect of a single term in the anisotropic regularizer \eqref{Reg} for different values of $\sigma$ and $\tilt$. The horizontal and vertical components of the model gradient vary along the respective axes. Each ball represents the set of points where the regularization function equals a given value. Bottom row: Illustration of the behavior of the solution of a linearized version of \eqref{FWI_const} for $\tilt=0$ and various values of $\sigma$. Finding a solution 
is equivalent to finding a point constrained to a given line while minimizing the anisotropic regularization term; this is done by expanding the regularization ball until it becomes tangent to the line (see \cite[Figure 4.3]{Tarantola_2005_IPT}). 
}
 }
 \label{reg_balls}
 \vspace{-0.5cm}
 \end{figure}
 
We consider the following extended regularization functional:
\begin{equation} \label{Reg}
\mathcal{R}(m,\theta,\sigma)=\frac12\sum_{i=1}^n\|\Sigma_i \bold{R}([\bolds{\tilt}]_i)[\nabla \bold{m}]_i\|_2^2,
\end{equation}
where the tilt angles $\theta \in {\mathbb{R}}^{n \times 1}$ and weight matrices $\Sigma_i$ are optimally selected for each pixel of the model. The tilt angle $[\bolds{\tilt}]_i \in (-\frac{\pi}{2},\frac{\pi}{2}]$ represents the local structural orientation at pixel $i$, measured clockwise from the horizontal axis.
The rotation matrix $\bold{R}$ defined as:
\begin{equation}
\bold{R}([\bolds{\tilt}]_i)=
\begin{pmatrix}
~~\cos[\bolds{\tilt}]_i & \sin[\bolds{\tilt}]_i\\
-\sin[\bolds{\tilt}]_i & \cos[\bolds{\tilt}]_i
\end{pmatrix},
\end{equation}
which rotates the gradient components $[\nabla\bold{m}]_i$ \eqref{grad_comp} to align with the local structure's direction \tcm{and its normal}.
The weighting matrix $\Sigma_i$ is given by:
\begin{equation}
\Sigma_i=
\begin{pmatrix}
[\sigma_{x}]_i & 0\\
0 & [\sigma_{z}]_i
\end{pmatrix},
\end{equation}
where $[\sigma_{x}]_i, [\sigma_{z}]_i \geq 0$ and $[\sigma_{x}]_i+[\sigma_{z}]_i=1$.
These weights control the smoothing strength along the tilt direction $[\bolds{\tilt}]_i$ and its orthogonal direction. While \cite{Gholami_2024_RES} employed fixed values for these weights, which may lead to suboptimal results, we choose the weights that minimize the regularization function. 
Considering the constraint $[\sigma_{x}]_i+[\sigma_{z}]_i=1$, the weight matrix $\Sigma_i$ is determined by a single parameter $[\sigma]_i$ with $[\sigma_{x}]_i=[\sigma]_i$ and $[\sigma_{z}]_i=1-[\sigma]_i$.

Using the regularization \eqref{Reg} the FWI problem becomes
\begin{mini}
{m,u_s,\theta,\sigma}{\frac{1}{2}\sum_{s=1}^{n_s}\|Pu_s-d_s\|_2^2+\alpha \mathcal{R}(m,\theta,\sigma)}{}{}
\addConstraint{\begin{array}{l}
A(m)u_s=b_s,\quad s =1,\ldots,n_s\vspace{0.1cm}\\
-\frac{\pi}{2}<\theta\leq \frac{\pi}{2}
\end{array}}{}{}
\label{FWI_const}
\end{mini}
where $\alpha >0$. 
The FWI problem in \eqref{FWI_const} ensures that the estimated model, while satisfying the data (and the wave equation), is smooth along the direction $\theta$. 

Figure \ref{reg_balls}\tcm{, top row,} 
\tcm{shows the }
anisotropic regularization balls  for only one term of the regularizer \eqref{Reg} and for different values of $[\sigma]_i=\sigma$ and $[\theta]_i=\theta$. We observe that for $\sigma = 0.5$, the anisotropic regularization becomes equivalent to the standard isotropic regularization. When the value of $\sigma$ increases, reaching the maximum considered value of 0.9, the degree of regularization applied in the direction $\theta$ and its normal 
is maximally different, resulting in a needle-shaped ellipse. This shape favors models with elongated features aligned with $\theta$, while allowing variations in the normal direction to it\tcm{, as illustrated in the bottom row}.

\section{Algorithm}
The augmented Lagrangian function of \eqref{FWI_const} is
\begin{align}
&\mathcal{L}(\theta,\sigma,u_s,m,\lambda,\nu)= 
    \frac{1}{2}\sum_{s=1}^{n_s}\|Pu_s-d_s\|_2^2+\alpha \mathcal{R}(m,\theta,\sigma) \nonumber\\
    & \quad+ \frac{\mu}{2}\sum_{s=1}^{n_s}\|A(m)u_s-b_s\|_2^2-\lambda_s^T(A(m)u_s-b_s) \nonumber\\
    &\quad +\frac{\tau}{2}\|\theta-z\|_2^2-\nu^T(\theta-z),\nonumber 
\end{align}
where $z(\theta)=\min(\max(\theta,-\frac{\pi}{2}),\frac{\pi}{2})$, $\lambda_s$ and $\nu$ are Lagrange multipliers associated to the constraints, and $\mu>0$ and $\tau>0$ are the penalty parameters. The regularization parameter $\alpha>0$ is given, being typically determined by the discrepancy principle. 

This problem can be efficiently solved using the alternating direction method of multipliers (ADMM, \cite{Boyd_2011_DOS}), resulting in the following iterative procedure starting from initial guess $m$ (typically encoding some prior information), 
(spatially invariant) $\sigma=0.5$
and initial multipliers $\lambda_s=\nu=0$: 
\begin{subequations} \label{alg}
\begin{align}
&\theta^+ = \arg\min_{\theta}~ \mathcal{L}(\theta,\sigma,u_s,m,\lambda,\nu) \label{alg_theta} \\
&\sigma^+ = \arg\min_{\sigma }~ \mathcal{L}(\theta^+,\sigma,u_s,m,\lambda,\nu)  \label{alg_sigma} \\
&u^+_s = \arg\min_u~  \mathcal{L}(\theta^+,\sigma^+,u,m,\lambda,\nu)  \label{alg_u}\\
&m^+ = \arg\min_m~  \mathcal{L}(\theta^+,\sigma^+,u^+_s,m,\lambda,\nu)  \label{alg_m}\\
&\lambda^+_s=\lambda_s - \mu(\bold{A}(m^+)\bold{u}_s^+-\bold{b}_s) \label{alg_eps}\\
&\nu^+=\nu -\tau (\theta^+-z^+). \label{alg_eps}
\end{align}
\end{subequations}
Here, the updated variables at each iteration are denoted by a superscript ``+".
Each step of the algorithm addresses a specific subproblem, as described more in details below. 

Subproblem \eqref{alg_theta} involves minimizing the regularization function $\mathcal{R}$, as defined in \eqref{Reg}, with respect to $\theta$, given the current values of $m$ and $\sigma$. This can can be efficiently achieved using a single iteration of the Gauss-Newton method. To ensure the stability of the $\theta$ update, it is essential to incorporate a smoothing term, as suggested by \cite{Gholami_2024_RES}.

Subproblem \eqref{alg_sigma} involves minimizing $\mathcal{R}$ with respect to $\sigma$, resulting in: 
\begin{equation} \label{sigma}
    [\sigma]_i = \frac{[g_{z'}]_i^2}{[g_{x'}]_i^2+[g_{z'}]_i^2},
\end{equation}
where $\left([g_{x'}]_i, [g_{z'}]_i\right)^T = \bold{R}([{\theta}]_i) [\nabla \bold{m}]_i$ represent the rotated gradient of the current model at the \(i\)th pixel. This expression for \eqref{sigma} is justified by considering that, when $[{\theta}]_i$ corresponds to the correct orientation angle, we typically expect  $|[g_{x'}]_i| \ll |[g_{z'}]_i|$. In this case, directly minimizing $[g_{x'}]_i^2 + [g_{z'}]_i^2$ would heavily penalize $[g_{z'}]_i$, which is undesirable because we aim to apply more smoothing along the direction defined by $[{\theta}]_i$. The weights $[\sigma]_i$ in \eqref{sigma} balance the contributions of gradient components in the weighted sum $[\sigma]_i^2[g_{x'}]_i^2 + (1-[\sigma]_i)^2[g_{z'}]_i^2$, ensuring that the magnitudes of both terms are approximately equal before penalization. This adaptive weighting mechanism effectively enhances smoothing along the structural direction $[{\theta}]_i$ while allowing variations across it, eventually enabling the recovery of dominating structures and details in the model $m$.  
To maintain stability during the $\sigma$ update, the computed weights are smoothed using a convolution with a $5\times 5$ averaging filter; alternatively, one could incorporate a smoothing term.

A comprehensive analysis of the methods for solving subproblems  \eqref{alg_u} and \eqref{alg_m} can be found in \cite{Gholami_2024_FWI}. 
 \begin{figure}
 \center
 \includegraphics[width=1\columnwidth]{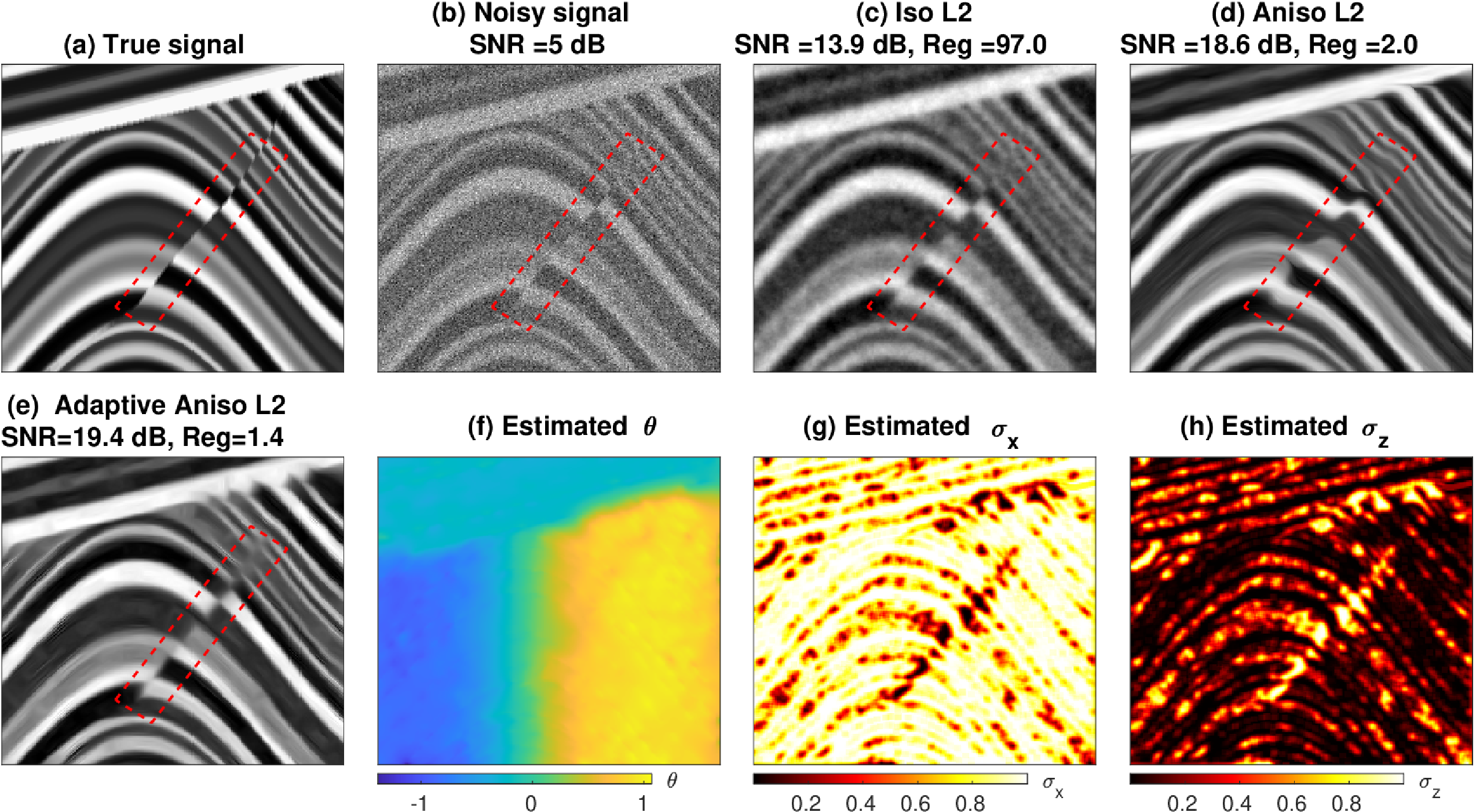}
 \caption{(a) True signal. (b) Noisy signal. Estimated signal by (c) the isotropic L2-norm, (d) the anisotropic L2-norm with fixed $\sigma$, (e) the adaptive anisotropic L2 with the estimated space-variying (f) $\theta$, (g) $\sigma_x$, and (h) $\sigma_z$.}
 \label{denoise}
 \vspace{-0.5cm}
 \end{figure}
\section{Numerical examples}
To evaluate the effectiveness of the proposed regularization approach, we conducted tests using first a denoising example and then FWI. 
\subsection{Denoising exmaple}
We considered the problem of random noise attenuation using the signal representing sedimentary layers, a plane unconformity, and a curved fault in Fig. \ref{denoise}(a). The noisy version of this signal, with a signal-to-noise ratio (SNR) of 5 dB, is shown in Fig. \ref{denoise}(b). The denoising task involves solving the optimization problem:
\begin{equation}\label{noise}
\minimize_{m}~ \frac{1}{2}\|m-d\|_{2}^{2}+\frac{\alpha}{2}\mathcal{R}(m,\theta,\sigma),
\end{equation}
where $d$ is the noisy signal and $\mathcal{R}$ is the regularization term defined in \eqref{Reg}.
We addressed this problem using the ADMM framework outlined in \cite{Gholami_2024_RES}. To assess the impact of different weighting strategies for $\sigma$, we explored three distinct approaches. For each case, the regularization parameter $\alpha$ was determined using the discrepancy principle. Specifically, $\alpha$ was chosen such that $\|m(\alpha) - d\| = \|e\|$, where $e$ represents the true noise added to the signal. 
\begin{enumerate}
\item Constant weights $\sigma=0.5$.  This approach corresponds to the standard isotropic Tikhonov regularization. 
The denoised result is shown in Fig. \ref{denoise}(c). 
The associated SNR value and the regularization function are 13.9 dB and 97.5, respectively.

\item Constant weights $\sigma=0.9$. Here, the emphasis is on applying stronger smoothing along the structural orientation $\theta$. The resulting denoised signal is presented in Fig. \ref{denoise}(d), with an associated SNR of 18.6 dB and a regularization function value of 2. These results indicate a notable improvement in SNR and a significant reduction in the regularization term, suggesting that the estimated tilt field effectively aligns with the true signal orientation. However, the use of fixed weights leads to spatially invariant smoothing kernels, which causes oversmoothing in areas with discontinuities (highlighted by the rectangle in Fig. \ref{denoise}(a)-(e)). This limitation can be mitigated by employing spatially adaptive weights. 

\item Space varying weights. The weights were dynamically estimated according to \eqref{sigma}, minimizing the regularization function. The resulting denoised signal and the corresponding estimated $\theta$, $\sigma_x=\sigma$ and $\sigma_z=1-\sigma$  are shown in Fig. \ref{denoise}(e)-(h), respectively. This approach yields an SNR of 19.4 dB and a regularization function value of 1.4. 
\end{enumerate}
These results demonstrate that optimizing the weights reduced the regularization term while improving the SNR.

\subsection{Sparse FWI with regularization}
\begin{figure}[t!]
 \center
 \includegraphics[width=1\columnwidth,trim={0 0cm 0 0cm},clip]{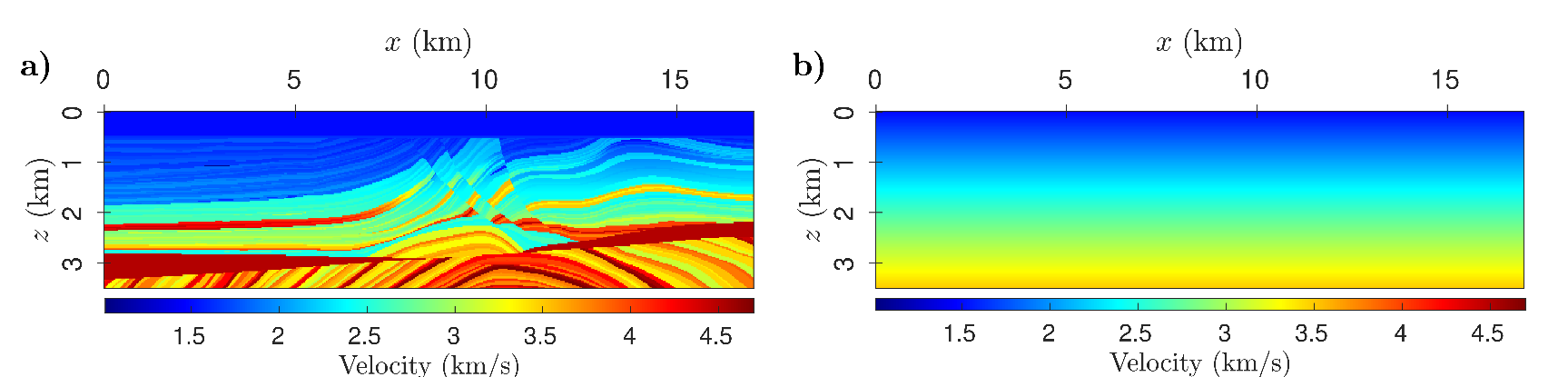}
 \caption{(a) The true Marmousi II velocity model; (b) The 1D starting velocity model.}
 \label{Marmousi_config}
 \vspace{-0.5cm}
 \end{figure}

For the FWI problem, we use the Marmousi II velocity model that is 17 km long and 3.5 km deep (Fig. \ref{Marmousi_config}a). The stationary-recording acquisition consists of 114 seismometers evenly spaced 150 meters apart along the seafloor, and 137 uniformly distributed pressure sources positioned at a depth of 25 meters. To improve computational efficiency, we apply the spatial reciprocity of Green's functions, treating sources as receivers and vice versa.

A perfectly matched layer (PML) absorbing boundary condition is implemented around the model to minimize boundary reflections, and a Ricker wavelet with a dominant frequency of 10 Hz is used as the source signature. The inversion starts with a 1D initial velocity model, which increases linearly with depth (Fig. \ref{Marmousi_config}b).
Inversion is performed in three cycles, following a classical multiscale frequency continuation approach, where frequencies progress from low to high. In each cycle, different frequencies are processed: from 3 to 6.5 Hz in the first cycle, 3 to 9 Hz in the second cycle, and 3 to 12 Hz in the third cycle, with a 0.5 Hz interval.
In each cycle, 10 iterations were performed to invert each frequency. The final model from the previous cycle serves as the initial model for the next cycle.

\begin{figure}
 \center
 \includegraphics[width=1\columnwidth]{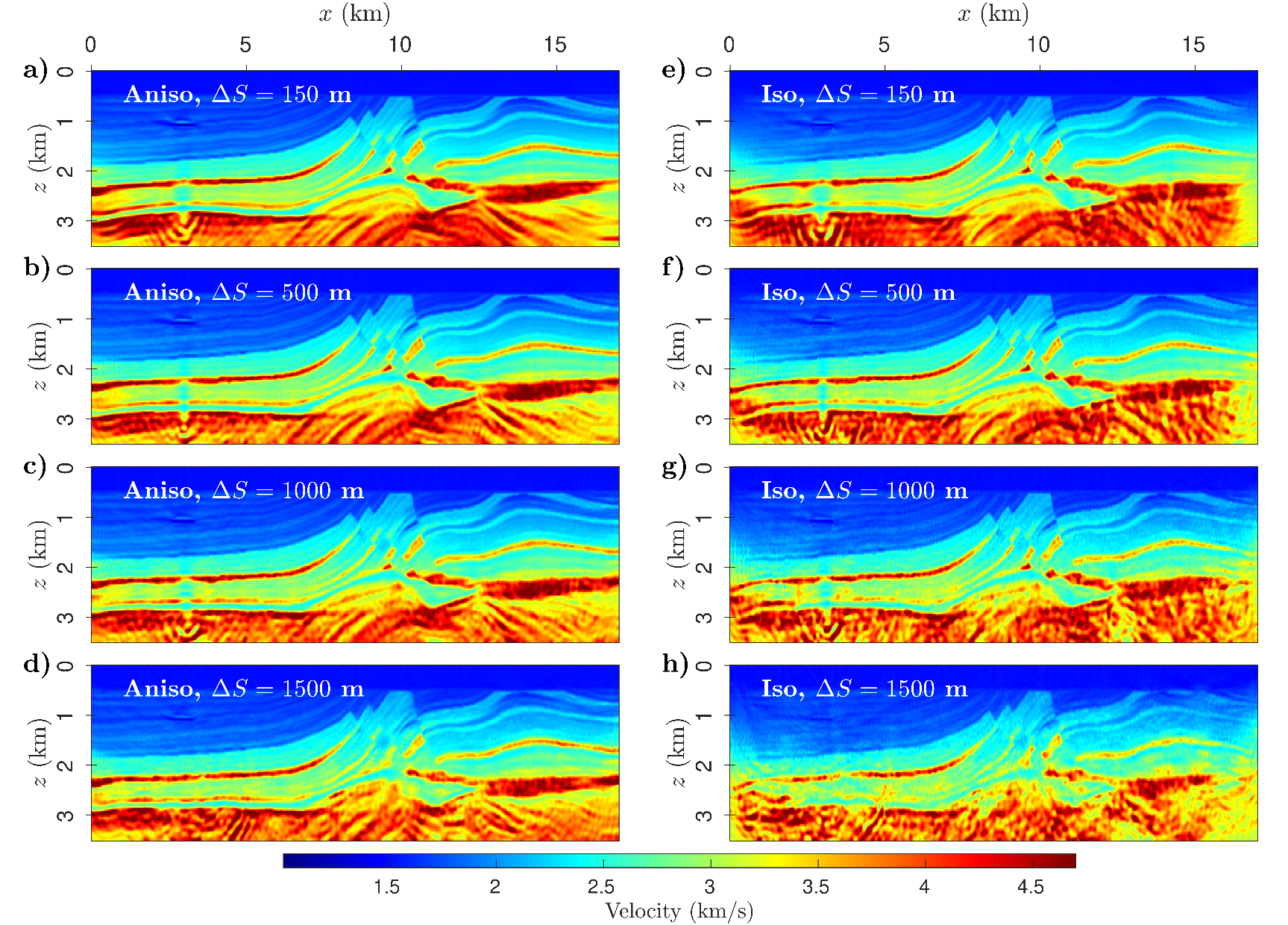}
 \caption{Comparison of FWI results using anisotropic and isotropic regularizations for different seismometer spacing ($\Delta S$). (a-d) anisotropic regularization for varying receiver spacings, (e-h) the corresponding isotropic results.}
 \label{Marm2_results}
 \vspace{-0.5cm}
 \end{figure}

We tested the performance of both anisotropic and isotropic regularizations in the inversion. The corresponding results are shown in Figs. \ref{Marm2_results}a and \ref{Marm2_results}e. The anisotropic regularization demonstrates superior performance, especially in the left, right, and bottom sections of the model, where subsurface illumination is weak.

To further highlight the role of regularization in sparse FWI, we increased the spacing between the ocean-bottom seismometers ($\Delta S$), transitioning to a sparse acquisition regime. We conducted tests with spacings of \tcm{500 m, 1000 m, and 1500 m, corresponding to 35, 18, and 11} receivers, respectively. The final inversion results are displayed in Figs. \ref{Marm2_results}b-d for anisotropic regularization and Figs. \ref{Marm2_results}f-h for isotropic regularization. The results show that anisotropic regularization is significantly more robust in handling sparse acquisition. 

 \begin{figure}
 \center
 \includegraphics[width=1\columnwidth]{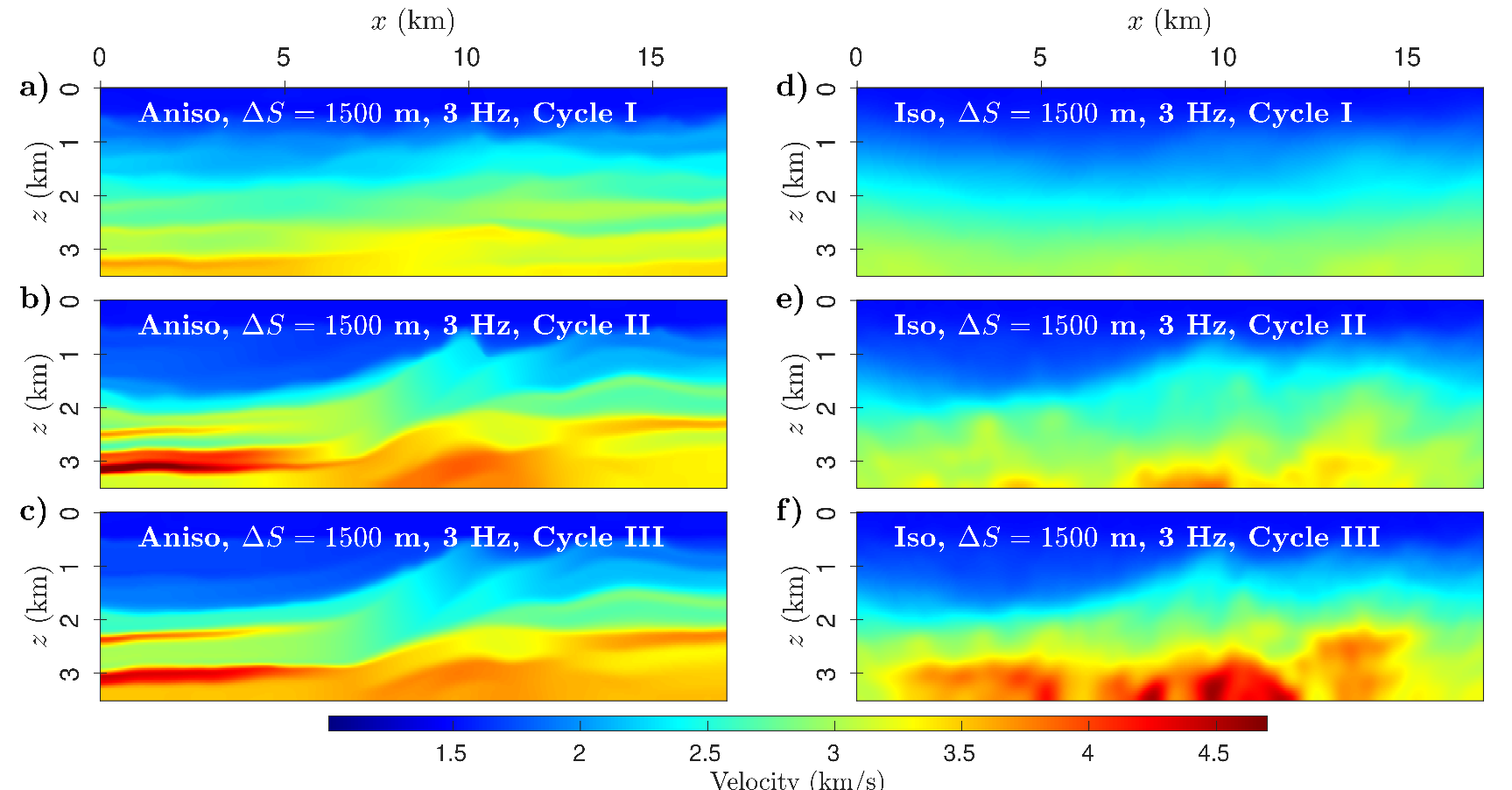}
 \caption{ Comparison of velocity models after inverting the first frequency (3 Hz) for the sparsest acquisition setup (\tcm{$\Delta S=1500$ m}). \tcm{(a)-(c)} Models obtained using anisotropic regularization for the first three inversion cycles. \tcm{(d)-(f)} Corresponding models obtained with standard isotropic regularization.}
 \label{Marm2_results_3hz}
 \end{figure}

For the sparsest acquisition case with \tcm{$\Delta S=1500$ m}, Figs. \ref{Marm2_results_3hz}\tcm{a-c} display the velocity models obtained after the anisotropic inversion of the first frequency (3 Hz) in each cycle. Figs. \ref{Marm2_results_3hz}\tcm{d-f} present the corresponding results using standard isotropic regularization. The adaptive, structure-oriented smoothing of anisotropic regularization leads to a significant improvement in model reconstruction, despite the very sparse sampling of the model.

\begin{figure}
 \center
 \includegraphics[width=0.75\columnwidth]{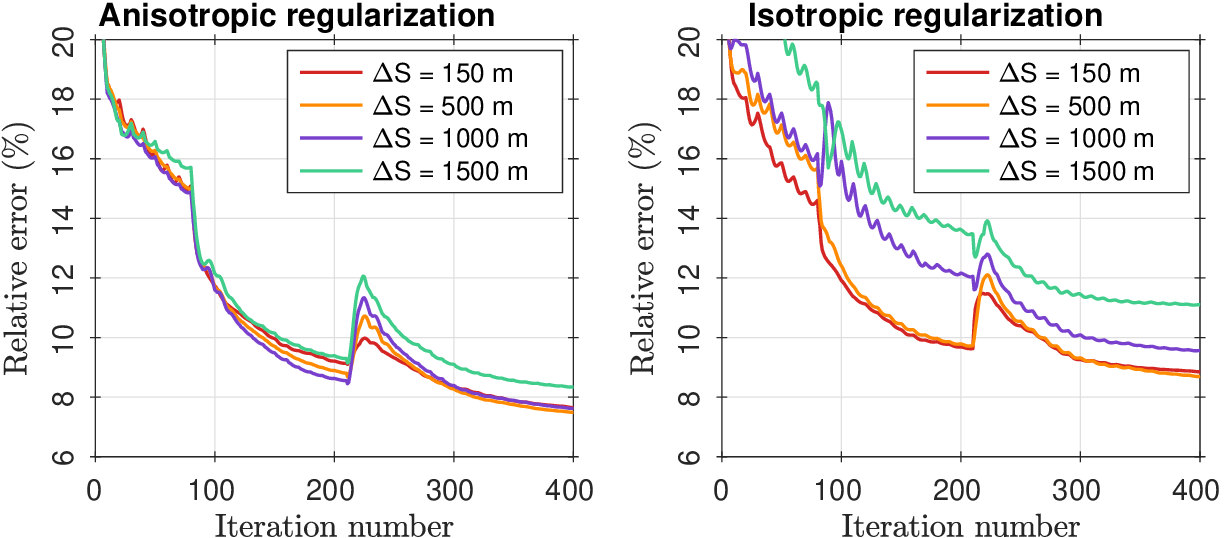}
 \caption{Evolution of the relative error during the ADMM FWI iterations for anisotropic and isotropic regularizations across different acquisition setups.}
 \label{Marm2_errors}
 \vspace{-0.5cm}
 \end{figure}

The superior performance of anisotropic regularization also is evident from the evolution of the model's relative error, defined as:
\begin{equation}
\text{Relative error}=100\times \frac{\|m_k-m^*\|_2}{\|m^*\|_2},
\end{equation}
where $m^*$ represents the true model, and $m_k$ is the model estimated at the $k$-th iteration. Figure \ref{Marm2_errors} illustrates the relative error for both anisotropic and isotropic regularizations, corresponding to the results shown in Figure \ref{Marm2_results}. \tcm{While the behavior of the error curve with respect to the number of sources is nonlinear, we observe that anisotropic regularization consistently outperforms isotropic regularization across all acquisition setups.} This difference is particularly pronounced in the sparsest acquisition setup with \tcm{$\Delta S = 1500$} m (the \tcm{green} curve in Figure \ref{Marm2_errors}). Anisotropic regularization produced a result comparable to those obtained with denser data, whereas isotropic regularization became trapped in a local minimum, as indicated by the error curve flattening between 300 and 400 iterations.
\tcm{Even if not reported here, we also observed that the performance of anisotropic regularization with a fixed weight falls between that of the locally optimized anisotropic regularization and the isotropic regularization.}

\section{Discussion and Conclusion}
Full waveform inversion (FWI) for estimating subsurface parameters is highly sensitive to various errors, including inaccuracies in the initial model, the absence of low-frequency content in the data, and sparse data acquisition leading to undersampling of the model. To address these challenges and stabilize the inversion, an appropriate regularization is essential, eventually ensuring a geologically meaningful solution. An effective regularization method should be able to identify the local orientations of structures within the model and apply adaptive smoothing along those directions, rather than across them.

We proposed an FWI method incorporating adaptive local anisotropic regularization. Unlike traditional methods, the orientation angles and the anisotropy weights in our approach are treated as part of the inverse problem and are estimated simultaneously with the unknown model parameters. This integration allows the algorithm to adapt dynamically to local structural orientations, leading to enhanced accuracy in complex geological environments.

Numerical experiments demonstrated that such anisotropic regularization proves especially effective for the inversion of sparse data, mitigating the impact of limited subsurface illumination. This capability underscores its potential for challenging practical applications, such as crustal-scale imaging using ocean-bottom seismometer (OBS) acquisition, where dense data acquisition is impractical or prohibitively expensive.
\section{Acknowledgments}  
This research was financially supported by the SONATA BIS grant
(No. 2022/46/E/ST10/00266) of the National Science Center in
Poland. 
%


\begin{thebibliography}{10}
\providecommand{\url}[1]{#1}
\csname url@samestyle\endcsname
\providecommand{\newblock}{\relax}
\providecommand{\bibinfo}[2]{#2}
\providecommand{\BIBentrySTDinterwordspacing}{\spaceskip=0pt\relax}
\providecommand{\BIBentryALTinterwordstretchfactor}{4}
\providecommand{\BIBentryALTinterwordspacing}{\spaceskip=\fontdimen2\font plus
\BIBentryALTinterwordstretchfactor\fontdimen3\font minus
  \fontdimen4\font\relax}
\providecommand{\BIBforeignlanguage}[2]{{%
\expandafter\ifx\csname l@#1\endcsname\relax
\typeout{** WARNING: IEEEtran.bst: No hyphenation pattern has been}%
\typeout{** loaded for the language `#1'. Using the pattern for}%
\typeout{** the default language instead.}%
\else
\language=\csname l@#1\endcsname
\fi
#2}}
\providecommand{\BIBdecl}{\relax}
\BIBdecl

\bibitem{Aghamiry_2019_CRO}
H.~Aghamiry, A.~Gholami, and S.~Operto, ``Compound regularization of
  full-waveform inversion for imaging piecewise media,'' \emph{{IEEE}
  Transactions on Geoscience and Remote Sensing}, vol.~58, no.~2, pp.
  1192--1204, 2020.

\bibitem{Jin_2022_RLM}
Y.~Jin, Y.~Zi, W.~Hu, Y.~Hu, X.~Wu, and J.~Chen, ``A robust learning method for
  low-frequency extrapolation in {GPR} full waveform inversion,'' \emph{IEEE
  Geoscience and Remote Sensing Letters}, vol.~19, pp. 1--5, 2022.

\bibitem{Operto_2023_FWI}
S.~Operto, A.~Gholami, H.~S. Aghamiry, G.~Guo, F.~Mamfoumbi, S.~Beller,
  K.~Aghazade, F.~Mamfoumbi, L.~Combe, and A.~Ribodetti, ``Extending the search
  space of full-waveform inversion beyond the single-scattering Born
  approximation: A tutorial review,'' \emph{Geophysics}, vol.~88, no.~6, pp.
  1--32, 2023.

\bibitem{Tian_2024_ESW}
W.~Tian, Y.~Liu, and X.~Di, ``Enhancing seismic waveform inversion using a
  three-step strategy with adversarial neural networks and seismic envelope,''
  \emph{IEEE Geoscience and Remote Sensing Letters}, 2024.

\bibitem{Metivier_2022_NAD}
L.~M{\'e}tivier and R.~Brossier, ``On the use of nonlinear anisotropic
  diffusion filters for seismic imaging using the full waveform,''
  \emph{Inverse Problems}, vol.~38, no.~11, p. 115001, 2022.

\bibitem{Barnier_2023_FWIp}
G.~Barnier, E.~Biondi, R.~G. Clapp, and B.~Biondi, ``Full waveform inversion by
  model extension: practical applications,'' \emph{Geophysics}, vol.~88, no.~5,
  pp. 1--138, 2023.

\bibitem{Wang_2020_ELW}
G.~Wang, T.~Alkhalifah, and S.~Wang, ``Enhancing low-wavenumber information in
  reflection waveform inversion by the energy norm born scattering,''
  \emph{IEEE Geoscience and Remote Sensing Letters}, vol.~19, pp. 1--5, 2020.

\bibitem{Wang_2019_RLW}
G.~Wang, S.~Yuan, and S.~Wang, ``Retrieving low-wavenumber information in FWI:
  An efficient solution for cycle skipping,'' \emph{IEEE Geoscience and Remote
  Sensing Letters}, vol.~16, no.~7, pp. 1125--1129, 2019.

\bibitem{Gorszczyk_2017_TRW}
A.~G{\'{o}}rszczyk, S.~Operto, and M.~Malinowski, ``Toward a robust workflow
  for deep crustal imaging by {FWI} of {OBS} data: The eastern nankai trough
  revisited,'' \emph{Journal of Geophysical Research: Solid Earth}, vol. 122,
  no.~6, pp. 4601--4630, Jun. 2017.

\bibitem{Tikhonov_1963_SIF}
A.~N. Tikhonov, ``Solution of incorrectly formulated problems and the
  regularization method.'' \emph{Sov Dok}, vol.~4, pp. 1035--1038, 1963.

\bibitem{Gholami_2024_RES}
A.~Gholami and S.~Gazzola, ``Robust estimation of structural orientation
  parameters and 2D/3D local anisotropic Tikhonov regularization,''
  \emph{Geophysics}, vol.~89, no.~6, pp. 1--81, 2024.

\bibitem{Pratt_1998_GNF}
R.~G. Pratt, C.~Shin, and G.~J. Hicks, ``{G}auss-{N}ewton and full {N}ewton
  methods in frequency-space seismic waveform inversion,'' \emph{Geophysical
  Journal International}, vol. 133, pp. 341--362, 1998.

\bibitem{Tarantola_2005_IPT}
A.~Tarantola, \emph{Inverse {P}roblem {T}heory and {M}ethods for {M}odel
  {P}arameter {E}stimation}.\hskip 1em plus 0.5em minus 0.4em\relax
  Philadelphia: Society for {I}ndustrial and {A}pplied {M}athematics, 2005.

\bibitem{Gholami_2024_FWI}
A.~Gholami and K.~Aghazade, ``Full waveform inversion and Lagrange
  multipliers,'' \emph{Geophysical Journal International}, vol. 238, no.~1, pp.
  109--131, 2024.

\bibitem{Boyd_2011_DOS}
S.~Boyd, N.~Parikh, E.~Chu, B.~Peleato, and J.~Eckstein, ``Distributed
  optimization and statistical learning via the alternating direction method of
  multipliers,'' \emph{Foundations and trends in machine learning}, vol.~3,
  no.~1, pp. 1--122, 2010.

\end{thebibliography}

\newcommand{\SortNoop}[1]{}

\end{document}